# Best rational approximations of an irrational number

Jean-Louis Sikorav

High Council for Economy
Ministry for the Economy and Finance
120 rue de Bercy
75772 Paris cedex 12, France

*Dedicated to the memory of Professor Denis Gerll 1913-2009*

*Abstract.* Given an irrational number $\alpha$ and an integer $q \neq 0$, there is a unique integer $p$, equal to the nearest integer to $q\alpha$, such that $|q\alpha - p| < \frac{1}{2}$. This result is used to derive theorems and algorithms for the best approximations of an irrational by rational numbers, which are improved with the pigeonhole principle and used to offer an informed presentation of the theory of continued fractions.

**Introduction.**

The approximation of a real number by a rational one is an ancient problem encountered in various branches of knowledge, illustrated in astronomy by the theory of calendars and in engineering by the design of cogwheel astronomical clocks. The problem also arose in arithmetic with the discovery of irrational numbers, such as the square root of $n$ when $n$ is a positive integer that is not an exact square. This finding has given birth to the field of rational - also called Diophantine - approximations of such numbers.

We know today that the set of rational numbers is - by construction - dense in the set of real numbers, so that for any given irrational number $\alpha$, and for an arbitrary small number $\varepsilon > 0$, there exist infinitely many fractions $\frac{p}{q}$, where $p$ and $q$ are integers with $q \neq 0$, that approximate $\alpha$ with the degree of accuracy $\varepsilon$: $\left|\alpha - \frac{p}{q}\right| < \varepsilon$. Among these fractions, best rational approximations of an irrational number, abbreviated here as *BRAIN*, are the solutions of specific problems of extrema:

*Given an irrational number $\alpha$ and a positive integer $n$, find the integers $(x, y)$, with $0 < y < n$, minimizing either $\left| \alpha - \frac{x}{y} \right|$ or $|y\alpha - x|$.*

The corresponding fractions $\frac{x}{y}$ are the best rational approximations of $\alpha$ having a denominator smaller in absolute value than a given quantity, known - using the terminology of Khinchin - as best approximations of the first kind (for the minimization of $\left| \alpha - \frac{x}{y} \right|$) and of the second kind (for the minimization of $|y\alpha - x|$).

The first complete solution of these two problems is associated with the works of Wallis, Huygens, Euler and Lagrange and a few others. It relies on the apparatus of continued fractions, itself an offspring of Euclid's algorithm, and provides existence theorems and algorithms for the computation of the approximations.

The present investigation contains:

1) A solution of the problem of best rational approximations using a straightforward approach based on the concept of nearest integer, also yielding existence theorems and algorithms;

2) An improvement of the knowledge gained, through a confrontation with Dirichlet's approximation theorem, leading to other theorems and algorithms;

3) An introduction to the field of continued fractions exploiting the results obtained in 1) and 2).

**1. Nearest integer of an irrational number**

- Let $\alpha$ be a positive irrational number. The set of natural integers greater than $\alpha$ is not empty: according to the well-ordering principle, here equivalent to a principle of descent, it possesses a least element. This defines the least integer greater than $\alpha$, denoted $\lceil \alpha \rceil$. The number denoted $\lfloor \alpha \rfloor = \lceil \alpha \rceil - 1$ is in turn the greatest integer less than $\alpha$. Both symbols have been conceived by Iverson. The two numbers $\lfloor \alpha \rfloor$ and $\lceil \alpha \rceil$ are the only integers $n$ such that $|\alpha - n| < 1$.

If $\alpha$ is a negative irrational number, $-\alpha$ is positive. We can write: $\lfloor -\alpha \rfloor < -\alpha < \lceil -\alpha \rceil$, hence: $-\lceil -\alpha \rceil < \alpha < -\lfloor -\alpha \rfloor$, showing in that case too the existence of a greatest integer less than $\alpha$ and a least integer greater than $\alpha$, with $\lfloor \alpha \rfloor = -\lceil -\alpha \rceil$ and $\lceil \alpha \rceil = -\lfloor -\alpha \rfloor$. This establishes the existence of such integers for all irrational numbers.

- We define the fractional part $\{\alpha\} = \alpha - \lfloor \alpha \rfloor$, with $0 < \{\alpha\} < 1$. This number measures the distance between $\alpha$ and $\lfloor \alpha \rfloor$ and the number $1 - \{\alpha\}$ the distance between $\alpha$ and $\lceil \alpha \rceil$. Of these two numbers, only one is strictly smaller than $\frac{1}{2}$, as their sum is equal to $1$ and the equality $1 - \{\alpha\} = \{\alpha\}$ is ruled out by the irrationality of $\alpha$. Accordingly, of the two numbers $\lfloor \alpha \rfloor$ and $\lceil \alpha \rceil$, only one, called the nearest integer to $\alpha$ and denoted here by the more symmetrical symbol $[\alpha]$ (Gauss square bracket) is such that:

$$0 < |\alpha - [\alpha]| < \frac{1}{2}.$$

The nearest integer verifies: $[\alpha] = \left\lfloor \alpha + \frac{1}{2} \right\rfloor = \left\lceil \alpha - \frac{1}{2} \right\rceil$.

- The distance between $\alpha$ and the nearest integer is written $\|\alpha\|$, with:

$$\|\alpha\| = |\alpha - [\alpha]| = \min(\{\alpha\}, 1 - \{\alpha\}).$$

## 2. Nearest rational of an irrational number, irrationality criteria and representations of irrational numbers

Let $q$ be a strictly positive integer. The preceding results can be generalized through the replacement of $\alpha$ by the irrational number $q\alpha$:

- There are exactly two integers $p$, equal to $\lfloor q\alpha \rfloor$ and $\lceil q\alpha \rceil$, such that:

$$|q\alpha - p| < 1.$$

For the first integer $|q\alpha - \lfloor q\alpha \rfloor| = \{q\alpha\}$ and for the second $|q\alpha - \lceil q\alpha \rceil| = 1 - \{q\alpha\}$.

- There is a unique integer $r = [q\alpha]$ such that:

$$0 < |q\alpha - r| = \|q\alpha\| < \frac{1}{2}.$$

This leads to basic irrationality criteria:

Let $\alpha$ be a real number. The following statements are equivalent:

(i)     $\alpha$ is irrational;

(ii)    For any strictly positive integer $q$, there are exactly two integers $p$, such that:

$$|q\alpha - p| < 1;$$

(iii)   For any strictly positive integer $q$, there is a unique integer $r$ such that:

$$0 < |q\alpha - r| = \|q\alpha\| < \tfrac{1}{2}.$$

Dividing by $q$ the members of the inequalities above further gives:

- There are exactly two integers $p$, such that:

$$\left|\alpha - \frac{p}{q}\right| < \frac{1}{q}.$$

- There is a unique integer $r$ such that:

$$0 < \left|\alpha - \frac{r}{q}\right| = \frac{\|q\alpha\|}{q} < \frac{1}{2q}.$$

These results lead to further irrationality criteria:

Let $\alpha$ be a real number. The following statements are equivalent:

(i)     $\alpha$ is irrational;

(ii)    For any strictly positive integer $q$, there are exactly two integers $p$ such that:

$$\left|\alpha - \frac{p}{q}\right| < \frac{1}{q};$$

(iii)   For any strictly positive integer $q$, there is a unique integer $r$ such that:

$$0 < \left|\alpha - \frac{r}{q}\right| < \frac{1}{2q}.$$

In this last inequality, one cannot lower the denominator $2q$ (say by replacing the factor $2$ by $(2 - \varepsilon)$ or $q$ by $q^{1-\epsilon}$, with $0 < \varepsilon < 1$) without losing the uniqueness of $r$.

The sequence $(\frac{[q\alpha]}{q})_{q\geq 1}$ is such that the distance between its $q$th term and $\alpha$ is smaller $\frac{1}{q}$; thus it converges towards $\alpha$ when $q$ goes to infinity: $\lim_{q\to\infty}\frac{[q\alpha]}{q} = \alpha$. The same conclusions apply to the sequences of underestimates $(\frac{\lfloor q\alpha\rfloor}{q})_{q\geq 1}$ and of overestimates $(\frac{\lceil q\alpha\rceil}{q})_{q\geq 1}$ of $\alpha$. These three sequences are two by two distinct and each provides a unique representation of $\alpha$.

The rational number $\frac{[q\alpha]}{q}$ is the nearest rational fraction to $\alpha$ for an arbitrary numerator and the fixed denominator $q$. Stated otherwise, it is the best rational approximation of $\alpha$ for this denominator. To acknowledge the ease with which this best approximation is derived, we shall call it $TRAIN\,(\alpha, q)$ (for $\underline{T}$rivial $\underline{A}$pproximation of an $\underline{I}$rrational $\underline{N}$umber). It will now be used to obtain best approximation theorems and algorithms.

## 3. *TRAIN* to *BRAIN*

### 3.1 *BRAIN Theorems*

*First BRAIN theorem.* Let $\alpha$ be an irrational number. There is a unique, infinite sequence of rationals $(\frac{p_k}{q_k})_{k\geq 1}$, where $p_k$ and $q_k$ are relatively prime integers with $1 \leq q_k$ and $p_k = [q_k\alpha]$, such that $\left|\alpha - \frac{p_k}{q_k}\right| = \min_{1\leq q\leq q_k}\frac{\|q\alpha\|}{q}$ and $\left|\alpha - \frac{p_{k+1}}{q_{k+1}}\right| < \left|\alpha - \frac{p_k}{q_k}\right|$. The fraction $\frac{p_k}{q_k}$ is the $k^{th}$ best rational approximation of the first kind of $\alpha$ and denoted $BRAIN\,I(\alpha, k)$. The sequence $(BRAIN\,I(\alpha, k))_{k\geq 1}$ converges to $\alpha$ when $k \to \infty$ and provides a unique representation of $\alpha$.

*Proof.*

Let $N$ be a strictly positive integer. For each integer $q$ with $1 \leq q \leq N$, the smallest value of $\left|\alpha - \frac{p}{q}\right|$ where $p$ is an integer, is equal to $\frac{\|q\alpha\|}{q}$ and is obtained for $p = [q\alpha]$: the fraction $\frac{p}{q}$ is simply $TRAIN\,(\alpha, q)$. This implies that the set of numbers of the form $\left|\alpha - \frac{p}{q}\right|$, where $p$ and $q$ are integers with $1 \leq q \leq N$ possesses a smallest element equal to the smallest of the $N$ numbers

$\frac{\|q\alpha\|}{q}$ with $1 \leq q \leq N$. This number, $\min_{1 \leq q \leq N} \frac{\|q\alpha\|}{q}$, is obtained for one (or several) integer(s) $q_i$. Let us call $q_0$ the smallest denominator. If the two numbers $[q_0 \alpha]$ and $q_0$ are not relatively prime, then we can write $[q_0 \alpha] = us$ and $q_0 = ut$, where $s, t$ and $u$ are integers with $u \geq 2$, showing that $\left|\alpha - \frac{s}{t}\right|$ is also equal to $\min_{1 \leq q \leq N} \frac{\|q\alpha\|}{q}$, with $t < q_0$, in contradiction with the definition of $q_0$. The two numbers $[q_0 \alpha]$ and $q_0$ are thus relatively prime. [1] The fraction $\frac{[q_0 \alpha]}{q_0}$ is therefore irreducible and is the best rational approximation of the first kind with a denominator $\leq N$. For $N = 1$, the first term of the sequence $BRAIN\ I(\alpha, k)$, $BRAIN\ I(\alpha, 1)$ is given by $q_1 = 1$ and $p_1 = [\alpha]$. Since the sequence $(\frac{[q\alpha]}{q})_{q \geq 1}$ converges towards $\alpha$ when $q \to \infty$, there is an integer $n' > N$ such that $\left|\alpha - \frac{[n'\alpha]}{n'}\right| < \left|\alpha - \frac{[q_0\alpha]}{q_0}\right|$. This shows that the sequence $(BRAIN\ I(\alpha, k))_{k \geq 1}$ is unbounded and converges to $\alpha$ when $k \to \infty$. Lastly, one can see that by construction, $k' > k$ implies $q_{k'} > q_k$ and $\left|\alpha - \frac{p_{k'}}{q_{k'}}\right| < \left|\alpha - \frac{p_k}{q_k}\right|$: the sequence $\min_{1 \leq q \leq q_k} \frac{\|q\alpha\|}{q}$ strictly decreases to zero when $q_k \to \infty$.

*Second BRAIN theorem.* Let $\alpha$ be an irrational number. There is a unique, infinite sequence of rationals $(\frac{p_k}{q_k})_{k \geq 1}$, where $p_k$ and $q_k$ are relatively prime integers with $1 \leq q_k$ and $p_k = [q_k \alpha]$, such that: $|q_k \alpha - p_k| = \min_{1 \leq q \leq q_k} \|q\alpha\|$ and $|q_{k+1}\alpha - p_{k+1}| < |q_k \alpha - p_k|$. The fraction $\frac{p_k}{q_k}$ is the $k^{th}$ best rational approximation of the second kind of $\alpha$ and denoted $BRAIN\ II(\alpha, k)$. The

---

[1] One can also reach this conclusion by noting that the equation $\left|\alpha - \frac{p}{q}\right| = \left|\alpha - \frac{r}{s}\right|$, where $\alpha$ is irrational and $p, q, r, s$ are integers (with $q > 0$ and $s > 0$), implies $\alpha - \frac{p}{q} = \alpha - \frac{r}{s}$, the other possibility $\alpha - \frac{p}{q} = -(\alpha - \frac{r}{s})$ being ruled out by the irrationality of $\alpha$. Thus $ps = qr$, and by choosing $p$ and $q$ relatively prime, we deduce that $p$ divides $r$ and $q$ divides $s$: there exists an integer $m \geq 1$ such that $r = mp$ and $s = mq$ (implying $s \geq q$).

sequence $(BRAIN\ II(\alpha, k))_{k \geq 1}$ converges to $\alpha$ when $k \to \infty$ and provides a unique representation of $\alpha$.

*Proof.*

Let $N$ be a strictly positive integer. For each integer $q$ with $1 \leq q \leq N$, the smallest value of $|q\alpha - p|$ where $p$ is an integer is equal to $\|q\alpha\|$ and is obtained for $p = [q\alpha]$. This implies that the set of numbers of the form $|q\alpha - p|$, where $p$ and $q$ are integers with $1 \leq q \leq N$ possesses a smallest element equal to the smallest of the $n$ numbers $\|q\alpha\|$ with $1 \leq q \leq N$. Given the irrationality of $\alpha$, it is easy to see that the equation $|q\alpha - p| = |s\alpha - r|$ where $p, q, r, s$ are integers (with $q > 0$ and $s > 0$) implies that $q = s$ and $p = r$ ($q \neq q' \Rightarrow \|q\alpha\| \neq \|q'\alpha\|$). The number $\min_{1 \leq q \leq N} \|q\alpha\|$ is therefore obtained for a unique integer $q_0$. If the two numbers $[q_0\alpha]$ and $q_0$ are not relatively prime, then we can write $[q_0\alpha] = us$ and $q_0 = ut$, where $s, t$ and $u$ are integers with $u \geq 2$, showing that $|t\alpha - s| = \frac{1}{u}|q_0\alpha - p_0|$, in contradiction with the definition of $q_0$. The two numbers $[q_0\alpha]$ and $q_0$ are thus relatively prime. The fraction $\frac{[q_0\alpha]}{q_0}$ corresponds to $TRAIN\ (\alpha, q_0)$. It is irreducible and is the best rational approximation of the second kind with a denominator $\leq N$.

Again by construction, $k' > k$ implies $q_{k'} > q_k$ and $|q_{k'}\alpha - p_{k'}| < |q_k\alpha - p_k|$. The proof that the sequence $(BRAIN\ II(\alpha, k))_{k \geq 1}$ is unbounded and converges to $\alpha$ is similar to that given above.

*Comparison of the two sequences $(BRAIN\ I(\alpha, k))_{k \geq 1}$ and $(BRAIN\ II(\alpha, k))_{k \geq 1}$*

Let $\alpha$ be an irrational number, $p$ and $q$ integers with $1 \leq q$. If $p$ differs from $[q\alpha]$, $\frac{p}{q}$ cannot be best rational approximation of $\alpha$, either of the first or the second kind. Supposing now that $p = [q\alpha]$, if the fraction $\frac{[q\alpha]}{q}$ is not a best rational approximation the first kind of order $q$ of $\alpha$, there is a fraction $\frac{[r\alpha]}{r}$, with $1 \le r < q$ such that $\left|\alpha - \frac{[r\alpha]}{r}\right| \le \left|\alpha - \frac{[q\alpha]}{q}\right|$, from which we can deduce:

$$|r\alpha - [r\alpha]| \le \frac{r}{q}|q\alpha - [q\alpha]| < |q\alpha - [q\alpha]|,$$

showing that $\frac{[q\alpha]}{q}$ cannot be a best rational approximation of the second kind of $\alpha$. By contraposition, we conclude that an approximation of an irrational number $\alpha$ of the second kind is also an approximation of the first kind.

### 3.2 BRAIN algorithms

These two theorems establish the existence and the uniqueness of best rational approximations of order $n$ of an irrational number $\alpha$ (of the first or the second kind). Taken in isolation they remain ineffective, pure existence theorems. Yet they can serve for the elaboration of algorithms provided that an adequate approximation of $\alpha$ is available through another approach. In such circumstances these theorems possess companion algorithms:

1) Let $N$ be a strictly positive integer with. The following algorithm determines $BRAIN\ I(\alpha, k_0)$, where $k_0$ is the largest integer such that $q_{k_0} \le N$ among the irreducible fractions $\frac{[q_k\alpha]}{q_k}$ minimizing $\frac{\|q_k\alpha\|}{q_k}$ with $1 \le q_k \le N$:

- For each integer $q$, $1 \le q \le N$, compute $[q\alpha]$ and $\frac{\|q\alpha\|}{q}$. Sort the $N$ terms $\frac{\|q\alpha\|}{q}$ in ascending order and among the fraction(s) minimizing $\frac{\|q\alpha\|}{q}$, retain the one having the smallest denominator $q_{k_0}$: $BRAIN\ I(\alpha, k_0) = \frac{[q_{k_0}\alpha]}{q_{k_0}}$.

The determination of $\frac{\|N\alpha\|}{N}$ requires the knowledge of an approximation of $\alpha$ with an accuracy of about $\frac{1}{N}$, and therefore with a number of digits proportional to $logN$. This in turn implies that computation of the product $N\alpha$ grows with the square of $logN$ and that the time needed to run the complete algorithm is proportional to $N(logN)^2$.

2) In the same manner, the following algorithm determines $BRAIN\ II(\alpha, k_0)$, where $k_0$ is the largest integer such that $q_{k_0} \leq N$ among the fractions $\frac{[q_k\alpha]}{q_k}$ minimizing $\|q_k\alpha\|$ with

$1 \leq q_k \leq N$:

- For each integer $q$, $1 \leq q \leq N$, compute $[q\alpha]$ and $\|q\alpha\|$. Sort the $N$ terms $\|q\alpha\|$ in ascending order, and retain the unique fraction minimizing $\|q\alpha\|$ characterized by its denominator $q_{k_0}$:

$BRAIN\ II(\alpha, k_0) = \frac{[q_{k_0}\alpha]}{q_{k_0}}$. The time needed to run the algorithm is proportional to $N(logN)^2$.

The two algorithms are easily implemented with a standard computer spreadsheet, such as Microsoft Excel, used here for convenience. We have illustrated them with the number $\pi$, equal to the ratio of the circumference to the diameter of a circle. For the sake of conciseness we omit the proof - due to Lambert and greatly simplified by Niven - of the irrationality of $\pi$.

Many dedicated algorithms have been described for the computation of approximations of $\pi$. The approach of Archimedes, based on regular inscribed and circumscribed polygons, yields the early classic results: $\frac{223}{71} < \pi < \frac{22}{7}$ (Archimedes), identifying $\frac{22}{7}$ as a fair approximation and also $3.1415926 < \pi < 3.1415927$ (Zu Chongzhi), leading to the identification of $\frac{355}{113}$ as a more accurate approximation.

The algorithm for $BRAIN\ I(\pi, k_0)$ has been run for the first $1000$ integers, with the twenty best results displayed in table 1. The smallest value for $\frac{\|q\pi\|}{q}$, $1 \leq q \leq 1000$, approximately equal to $2.66764 \times 10^{-7}$, is obtained with eight fractions, of which only one is irreducible: $\frac{355}{113}$, having the smallest denominator $113$, is the searched best approximation of the first kind. The seven other fractions are reducible to $\frac{355}{113}$ as predicted by the theorem and can be crossed out. The algorithm resembles the sieve of Eratosthenes used for the determination of prime numbers. Indeed, reducible fractions are erased here in the same manner as composite integers are excluded from the sieve of Eratosthenes. Only one best approximation can be determined in the twenty best results. A search for smaller denominators ($q < 113$) requires a complete inspection of the spreadsheet from the top to the bottom (not shown) leading to the determination of all the $BRAIN\ I(\pi, k_i)$ such that $q_{k_i} \leq 1000$, which, written in ascending order, reads:

$$\frac{3^+}{1}, \frac{13^-}{4}, \frac{16^-}{5}, \frac{19^-}{6}, \frac{22^-}{7}, \frac{179^+}{57}, \frac{201^+}{64}, \frac{223^+}{71}, \frac{245^+}{78}, \frac{267^+}{85}, \frac{289^+}{92}, \frac{311^+}{99}, \frac{333^+}{106}, \frac{355^-}{113}.$$

The symbols plus and minus are used to indicate underestimates ($(\pi - \frac{[q\pi]}{q}) > 0$) and overestimates ($(\pi - \frac{[q\pi]}{q}) < 0$) of $\pi$ respectively. The fraction $\frac{355}{113}$ is thus the 14$^{th}$ best approximation of the first kind ($\frac{355}{113} = BRAIN\ I(\pi, 14)$). We observe that $\frac{22}{7}, \frac{223}{71}$ (Archimedes' result) and $\frac{355}{113}$ (Zu Chongzhi's result) are best approximations of the first kind.

In the same manner, the algorithm for $BRAIN\ II(\pi, k_0)$ has been run for the first $1000$ integers, with the twenty best results being displayed in table 2. There is a unique fraction minimizing $\|q\pi\|, 1 \leq q \leq 1000$, as predicted by the theorem. It is equal to $\frac{355}{113}$, which is the searched best approximation of the second kind, with $\|355\pi\|$ being approximately equal to $3.0144 \times 10^{-5}$. The seven following best results are again fractions reducible to $\frac{355}{113}$. An inspection of the twenty best results leads now to an almost complete determination of all the $BRAIN\ II(\pi, k_0)$: only the first - $\frac{3}{1}$ - is buried deeper in the spreadsheet. The list of the $BRAIN\ II(\pi, k_0)$, written in ascending order, reads: $\frac{3^+}{1}, \frac{22^-}{7}, \frac{333^+}{106}, \frac{355^-}{113}$. We now identify $\frac{22}{7}$ and $\frac{355}{113}$ as the second and fourth best approximation of the second kind respectively.

## 4. Hidden *TRAIN*: Confrontation with Dirichlet's approximation theorem.

For an irrational number $\alpha$ and an integer $q \neq 0$, there is a unique solution $r = [q\alpha]$ to the inequality $\left|\alpha - \frac{r}{q}\right| < \frac{1}{2q}$ where $r$ is an integer. This result applies *a fortiori* if $\alpha$ satisfies a relation $\left|\alpha - \frac{r}{q}\right| < f(q)$, where $r$ and $q \neq 0$ are two integers and $f$ is a function of $q$ such that $f(q) < \frac{1}{2q}$, a situation frequently found in the field of Diophantine approximation. A paradigmatic case is Dirichlet's approximation theorem:

*Let $\alpha$ be an irrational number: there exist infinitely rational numbers $\frac{p}{q}$ such that: $\left|\alpha - \frac{p}{q}\right| < \frac{1}{q^2}$.*

The standard proof of this theorem relies on the Schubfachprinzip (literally principle of the drawers, principe des tiroirs in French, box or pigeonhole principle in English), which states that if $N + 1$ objects are put in $N$ drawers, at least one drawer will contain at least two objects.

*Proof*:

Let $N$ be a positive integer. The $N + 1$ numbers: $0, \{\alpha\}, \ldots \{k\alpha\}, \ldots \{(N-1)a\}, \{N\alpha\}$ are distributed over the $N$ disjoint intervals: $\left[0, \frac{1}{N}\right[ \ldots \left[\frac{i}{N}, \frac{i+1}{N}\right[ \ldots \left[\frac{N-1}{N}, 1\right[$, with $0 \leq i \leq N - 1$. The pigeonhole principle shows here that one interval, say $\left[\frac{i_0}{N}, \frac{i_0+1}{N}\right[$, must contain at least two of the $N + 1$ numbers, say $\{k\alpha\}$ and $\{l\alpha\}$. We can suppose without loss of generality $k > l$.

We have $|\{k\alpha\} - \{l\alpha\}| < \frac{1}{N}$ and since $\{k\alpha\} = k\alpha - \lfloor k\alpha \rfloor$ and $\{l\alpha\} = l\alpha - \lfloor l\alpha \rfloor$ we can write:

$|(k-l)\alpha - (\lfloor k\alpha \rfloor - \lfloor l\alpha \rfloor)| < \frac{1}{N}$, or, using $q = (k - l)$ and $p = \lfloor k\alpha \rfloor - \lfloor l\alpha \rfloor$: $|q\alpha - p| < \frac{1}{N}$ and also $\left|\alpha - \frac{p}{q}\right| < \frac{1}{qN}$ where $\frac{p}{q}$ is a rational with $1 \leq q \leq N$, which implies $\left|\alpha - \frac{p}{q}\right| < \frac{1}{q^2}$. This proves the existence of at least one rational number $\frac{p}{q}$, $1 \leq q \leq N$ for which the inequality $\left|\alpha - \frac{p}{q}\right| < \frac{1}{q^2}$ is satisfied.

Let us now suppose that there is only a finite number of fractions $\frac{p_i}{q_i}$ verifying the inequality above. Denoting $\varepsilon_i = \left| \alpha - \frac{p_i}{q_i} \right|$, we call $\varepsilon \, (> 0)$ the smallest of all $\varepsilon_i$. We choose $N > \frac{1}{\varepsilon}$: we know that there is a rational number $\frac{s}{t}$ such that $\left| \alpha - \frac{s}{t} \right| < \frac{1}{tN} \leq \frac{1}{N} < \varepsilon$. The fraction $\frac{s}{t}$ differs from all $\frac{p_i}{q_i}$, in contradiction with the hypothesis that their number is finite. This completes the proof.

As $\frac{1}{q^2}$ is smaller than $\frac{1}{2q}$ if $q \geq 2$, comparing this theorem with our previous results shows that for a fixed denominator $q \geq 2$, there is a unique numerator $p = [q\alpha]$ such that $\left| \alpha - \frac{p}{q} \right| < \frac{1}{q^2}$. We can thus state both a more precise version of Dirichlet's theorem and introduce a companion algorithm:

1) *Third BRAIN theorem.* Let $\alpha$ be an irrational number. There is a unique, infinite sequence of rationals $(\frac{p_k}{q_k})_{k \geq 1}$, where $p_k$ and $q_k$ are relatively prime integers with $1 \leq q_k$ and $p_k = [q_k\alpha]$, such that: $|q_k\alpha - p_k| = \|q_k\alpha\| < \frac{1}{q_k}$ and $k' > k$ implies $q_{k'} > q_k$. We shall call the fraction $\frac{p_k}{q_k}$ the $k^{th}$ best rational approximation of the third kind of $\alpha$, denoted $BRAIN\ III(\alpha, k)$. The sequence $(BRAIN\ III(\alpha, k))_{k \geq 1}$ converges to $\alpha$ when $k \to \infty$ and provides a unique representation of $\alpha$.

2) The following algorithm determines $BRAIN\ III(\alpha, k_0)$, where $k_0$ is the largest integer such that $q_{k_0} \leq N$ among the fractions $\frac{[q_k\alpha]}{q_k}$ such that $|q_k\alpha - p_k| < \frac{1}{q_k}$ with $1 \leq q_k \leq N$:

- For each integer $q$, $1 \leq q \leq N$, compute $[q\alpha]$ and $q\|q\alpha\|$. Sort the $N$ terms $q\|q\alpha\|$ in ascending order, and among those satisfying $q\|q\alpha\| < 1$, retain the unique fraction minimizing

$q\|q\alpha\|$ characterized by its denominator $q_{k_0}$: $BRAIN\ III(\alpha, k_0) = \frac{[q_{k_0}\alpha]}{q_{k_0}}$. The time needed to run the algorithm is again proportional to $N(logN)^2$.

The algorithm $BRAIN\ III(\pi, k_0)$ has been run for the first $\mathbf{1000}$ integers: there are 16 fractions such that $q\|q\pi\| < 1$, shown in table 3. The irreducible fraction $\frac{355}{113}$ is the searched best approximation of the third kind, with $355\|355\pi\|$ being approximately equal to $\mathbf{0.0034}$. The seven following best results are again fractions reducible to $\frac{355}{113}$. An inspection of the 16 results leads to a complete determination of all the $BRAIN\ III(\pi, k_0)$ in one sweep. The list of the $BRAIN\ III(\pi, k_0)$, written in ascending order, reads: $\frac{3^+}{1}, \frac{19^-}{6}, \frac{22^-}{7}, \frac{333^+}{106}, \frac{355^-}{113}$. It is mixture of all the best approximations of the second kind and a single best approximation of the second kind ($\frac{19^-}{6}$). This finding is explained by the two following theorems:

1) A best approximation of the second kind is a best approximation of the third kind.

*Proof*:

Let $N$ be a positive integer and $BRAIN\ II(\alpha, k_0) = \frac{[q_{k_0}\alpha]}{q_{k_0}}$ the best approximation of the second kind where $k_0$ is the largest integer such that $q_{k_0} \leq N$. From the demonstration of Dirichlet's theorem given above, there exist two integers $p$ and $q$, $1 \leq q \leq N$ such that $|q\alpha - p| < \frac{1}{N}$. By definition, $BRAIN\ II(\alpha, k_0)$ is such that $|q_{k_0}\alpha - [q_{k_0}\alpha]| < |q\alpha - p|$. It implies:

$|q_{k_0}\alpha - [q_{k_0}\alpha]| < \frac{1}{N}$ and *a fortiori* $|q_{k_0}\alpha - [q_{k_0}\alpha]| < \frac{1}{q_{k_0}}$; thus $\left|\alpha - \frac{[q_{k_0}\alpha]}{q_{k_0}}\right| < \frac{1}{q_{k_0}^2}$, showing that it is also a best approximation of the third kind. This explains why the list above contains all the best approximations of the second kind.

2) A best approximation of the third kind is a best approximation of the first kind.

*Proof:*

Let $N$ be a positive integer and $BRAIN\ III(\alpha, k_0) = \frac{[q_{k_0}\alpha]}{q_{k_0}}$ the best approximation of the third kind where $k_0$ is the largest integer such that $q_{k_0} \leq N$. If $\frac{[q_{k_0}\alpha]}{q_{k_0}}$ is not a best approximation of the first kind, there exits integers $r$ and $s$, $1 \leq s < q_{k_0}$ such that $\left|\alpha - \frac{r}{s}\right| < \left|\alpha - \frac{[q_{k_0}\alpha]}{q_{k_0}}\right| < \frac{1}{q_{k_0}^2}$. But this means that the distance between the two fractions $\frac{r}{s}$ and $\frac{[q_{k_0}\alpha]}{q_{k_0}}$ is smaller than $\frac{1}{q_{k_0}^2}$: $\left|\frac{r}{s} - \frac{[q_{k_0}\alpha]}{q_{k_0}}\right| < \frac{1}{q_{k_0}^2}$ or $\left|rq_{k_0} - s[q_{k_0}\alpha]\right| < \frac{s}{q_{k_0}} < 1$, which in turn implies that the nonzero integer $\left|rq_{k_0} - s[q_{k_0}\alpha]\right|$ is strictly smaller than 1. We conclude that $BRAIN\ III(\alpha, k_0)$ is also a best approximation of the first kind.

Beyond Dirichlet's theorem, the search for hidden TRAIN can be performed for other similar non-effective approximation theorems encountered the field of Diophantine approximation.

As a last comment, we note that the standard proof of Dirichlet's theorem uses only the $N$ fractional parts $\{k\alpha\}$, $1 \leq k \leq N$. The absence of the $N$ fractional parts $1 - \{k\alpha\}$ suggests the following improvement of the theorem:

*Let $\alpha$ be an irrational number: there exist infinitely rational numbers $\frac{p}{q}$ such that:* $\left|\alpha - \frac{p}{q}\right| < \frac{1}{2q^2}$

This theorem is actually a well-known result of the theory of continued fractions to which we turn our attention.

## 5. An informed presentation of the continued fraction apparatus

A common introduction of the theory of continued fractions starts with a definition of the regular continued fraction algorithm, followed by a sequence of deductions uncovering progressively its properties. The three BRAIN theorems and algorithms can be used to complete this standard introduction as described briefly in this section, which assumes a general knowledge of the basic theory of continued fractions and related issues, and where no attempt is made to prove the conclusions reached mostly empirically.

*Representation of an irrational number by a regular continuous fraction:*

Let $\alpha$ be an irrational number. We write first $\alpha = a_0 + \{\alpha\}$, with $a_0 = \lfloor \alpha \rfloor$, $0 < \{\alpha\} < 1$, then $a_1 = \left\lfloor \frac{1}{\{\alpha\}} \right\rfloor$ such that: $\frac{1}{\{\alpha\}} = a_1 + \left\{ \frac{1}{\{\alpha\}} \right\}$ and iterate the procedure: $a_2 = \left\lfloor \frac{1}{\left\{ \frac{1}{\{\alpha\}} \right\}} \right\rfloor$ and so forth. We can now write:

$$\alpha = a_0 + \{\alpha\} = a_0 + \cfrac{1}{a_1 + \left\{ \cfrac{1}{\{\alpha\}} \right\}} = a_0 + \cfrac{1}{a_1 + \cfrac{1}{a_2 + \cfrac{1}{a_3} \cdots}}$$

The irrationality of $\alpha$ implies that this process, called the regular continued fraction (RCF) algorithm, never ends. The integers $a_0, a_1, a_2, a_3, \ldots$ constitute the partial quotients of the continued fraction. They are strictly positive with the possible exception of $a_0$. The convergent of order $n$ of the continued fraction are the fractions $\frac{p_n}{q_n}$ defined by terminating the continued fraction with the partial quotient $a_n$: $\frac{p_0}{q_0} = \frac{a_0}{1}$ and for $n \geq 1$:

$$\frac{p_n}{q_n} = a_0 + \cfrac{1}{a_1 + \cfrac{1}{a_2 + \cfrac{1}{a_3 + \cdots \cdot \cfrac{1}{a_n}}}}$$

Let us use the number $\pi$ to investigate this algorithm. We find:

$$\pi = 3 + \cfrac{1}{7 + \cfrac{1}{15 + \cfrac{1}{1 + \cfrac{1}{292 + \cfrac{1}{1 + \cdots}}}}}$$

This gives the following fractions for the first six convergents:

$$\frac{p_0}{q_0} = \frac{3^+}{1}, \frac{p_1}{q_1} = \frac{22^-}{7}, \frac{p_2}{q_2} = \frac{333^+}{106}, \frac{p_3}{q_3} = \frac{355^-}{113}, \frac{p_4}{q_4} = \frac{103933^+}{33102}, \frac{p_5}{q_5} = \frac{104348^-}{33215}$$

A comparison with our previous results shows that the first four convergents are in fact the first four best approximation of the second kind of $\pi$. This suggests a first generalization that applies to subsequent convergents, such as $\frac{p_4}{q_4}$ and $\frac{p_5}{q_5}$, that we can check empirically using the *BRAIN* algorithm associated with the second *BRAIN* theorem. This confirms that:

$\frac{p_4}{q_4} = \frac{103933^+}{33102} = BRAIN\ II(\pi, 5)$, with $\|33102\pi\| \approx 1.91292 \times 10^{-5}$ and

$\frac{\|33102\pi\|}{33102\pi} \approx 5.77888 \times 10^{-10}$.

$\frac{p_5}{q_5} = \frac{104348^-}{33215} = BRAIN\ II(\pi, 6)$, with $\|33215\pi\| \approx 1.10151 \times 10^{-5}$ and

$\frac{\|33215\pi\|}{33215\pi} \approx 3.31631 \times 10^{-10}$.

We can therefore conjecture the basic approximation theorem:

*The sequence of the convergents $\left(\frac{p_k}{q_k}\right)_{k \geq 1}$ of an irrational number $\alpha$ is equal to the sequence $(BRAIN\ II(\alpha, k))_{k \geq 1}$.*

Two main facts emerge from a comparison of the *BRAIN* approach with the continued fraction theory:

1) Seemingly disjoint fractions $(BRAIN\ II(\alpha, k))_{k \geq 1}$ are now connected, being part of an infinite concatenation. An endless continued fraction is a mathematical oxymoron, perhaps best seen in the German name Kettenbruch (literally broken chain). The continued fraction algorithm is a discrete thread guiding us through the labyrinth of the continuum.

2) The efficiency of the continued fraction algorithm is greatly superior to that of the **BRAIN II** algorithm: the number of operations required to compute $\frac{p_4}{q_4}$ decreases from more than 33000 to just a few. Indeed, the algorithmic complexity is reduced by a factor $N$, from $N(log N)^2$ to $(log N)^2$.

The *BRAIN* approach further indicates that the choice of the floor function $\lfloor \alpha \rfloor$ in the regular continued fraction algorithm is somewhat arbitrary; the same conclusion holds for the use the ceiling function $\lceil \alpha \rceil$. A more reasonable choice would be the nearest integer function. When applied to the number $\pi$, the nearest integer continued fraction (NICF) algorithm yields:

$$\pi = 3 + \cfrac{1}{7 + \cfrac{1}{16 - \cfrac{1}{294 - \cfrac{1}{3 - \cfrac{1}{3 - \cdots}}}}}$$

This gives the following fractions for the first four convergents:

$$\frac{p_0}{q_0} = \frac{3}{1}, \frac{p_1}{q_1} = \frac{22}{7}, \frac{p_2}{q_2} = \frac{355}{113}, \frac{p_3}{q_3} = \frac{104348}{33215}$$

This last algorithm appears to speed up the regular continued fraction algorithm, at the expense of the completeness of type II best approximations.

Another conclusion drawn from the *BRAIN* approach is that in the convergents $\left(\frac{p_k}{q_k}\right)_{k \geq 1}$ of an irrational number $\alpha$, the numerator and the denominator are always related through the equation $p_k = [q_k \alpha]$. Applied to Archimedes' approximation of $\pi$, this gives:

$$\frac{p_1}{q_1} = \frac{22}{7} = \frac{[7\pi]}{7}$$

The fraction is obtained directly through the relationship between the numerator and the denominator. In comparison, with the RCF algorithm:

$$\frac{22}{7} = \lfloor \pi \rfloor + \left\lfloor \frac{1}{\{\pi\}} \right\rfloor = 3 + \frac{1}{7}.$$

One passes through intermediate stages in the computation of the reduced fraction. In the spirit of a problem of the mathematical Olympiads, this suggests a summary of the present work as a riddle:

*The fraction $\frac{22}{7}$ is a fair approximation of the number $\pi$, ratio of the circumference to the diameter of a circle. Find a simple relation between its numerator and denominator.*

Another way to learn more of the phenomenon of a continued fraction is to study its simplest realization, obtained through a saturation with compatible symmetries. Setting all the partial quotients $a_0, a_1, a_2, a_3$ ... equal to $1$ leads to the continued fraction:

$$1 + \cfrac{1}{1 + \cfrac{1}{1 + \cfrac{1}{1 + \cdots}}}$$

Looking into this fraction at greater magnifications, we observe the reappearance of the whole fraction itself (*Eadem resurgit*). Taking into account this property of scale invariance, one expects the unknown value $x$ of this unending fraction (of which the convergence is assumed) to satisfy the equation:

$$x = 1 + \frac{1}{x}$$

This quadratic equation has a single positive solution, $\varphi = \frac{1+\sqrt{5}}{2}$, known as the Golden ratio, towards which the continued fraction presumably converges.

We can take advantage of the *BRAIN* approach to investigate the best approximations of the Golden ratio. The algorithms for $BRAIN\ I(\varphi, k_0)$ and $BRAIN\ II(\varphi, k_0)$ have been run for the first $1000$ integers, with the twenty best results displayed in tables 4 and 5. There is a common best approximation for the two algorithms, equal to $\frac{1597^-}{687}$. An inspection of the full spreadsheets reveals actually that the 15 first two types of approximation are identical, being given in the following list:

$$\frac{2^-}{1}, \frac{3^+}{2}, \frac{5^-}{3}, \frac{8^+}{5}, \frac{13^-}{8}, \frac{21^+}{13}, \frac{34^-}{21}, \frac{55^+}{34}, \frac{89^-}{55}, \frac{144^+}{89}, \frac{233^-}{144}, \frac{377^+}{233}, \frac{610^-}{377}, \frac{987^+}{610}, \frac{1597^-}{987}$$

We are thus lead to hypothesize that this identity is always valid. We also notice that the approximations are successively over and underestimates. Looking back to the results obtained for $\pi$, one can see that this is also true for type II best approximations, but not for type I, suggesting that this result is of general validity for type II approximations of all irrational numbers.

Another conspicuous feature of the list given above is that the numerator of the $k^{th}$ approximation becomes the denominator of the $(k+1)^{th}$ approximation. Taking again into account the fact that every *BRAIN* is a *TRAIN* suggests that the list above is associated with the sequence defined recursively by:

$$\mathcal{F}_1 = 1 \text{ and } \mathcal{F}_{n+1} = [\mathcal{F}_n \varphi] \text{ for } n \geq 1:$$

$$1,2,3,5,8,13,21,34,55,89,144,233,377,610,987,1597\ldots$$

The recursive relation $\mathcal{F}_{n+1} = [\mathcal{F}_n \varphi]$ reveals the geometric character of the sequence. This is best seen in the equation:

$$\mathcal{F}_n = \left[\frac{\varphi^n}{\sqrt{5}}\right]$$

which is easily derived from the formula $\mathcal{F}_n = \frac{\varphi^n - (-\varphi)^{-n}}{\sqrt{5}}$.

When applied to the Golden ratio, the nearest integer continued fraction (NICF) algorithm yields: $\varphi = 2 - \cfrac{1}{3 - \cfrac{1}{3 - \cfrac{1}{3-1\ldots}}}$. The convergents of this continued fraction are $\frac{2^-}{1}, \frac{5^-}{3}, \frac{13^-}{8} \ldots$. They appear to skip every best approximation of even order, speeding up the rate of convergence by a factor 2.

The regular continued fraction algorithm offers another perspective on the sequence $\mathcal{F}_n$. Indeed, in the computation the convergents, the following recurrence relation shows up:

$$\mathcal{F}_{n+2} = \mathcal{F}_{n+1} + \mathcal{F}_n$$

The sequence $\mathcal{F}_n$ now appears as a sequence of integers defined by an additive recurrence relation, known as the Fibonacci sequence. The superiority of the RCF algorithm is again patent:

- In contrast with the *BRAIN* approach, which necessitates the independent knowledge of an approximation of the Golden ratio, the RCF algorithm is an "auto-algorithm", where an approximation of the Golden ratio is directly obtained.

- The first 15 best approximations are obtained must faster than with the *BRAIN* approach.

A merit of the *BRAIN* approach is that the Fibonacci numbers emerge without any computation of the additive recurrence relation.

We shall end our exploration of the properties of the best approximations of the Golden ratio with the results given by the $BRAIN\ III(\varphi, k_0)$, run for the first $1000$ integers and displayed in table 6. There are 15 fractions such that $q\|q\varphi\| < 1$. The fractions are in fact all such that $q\|q\varphi\| < \frac{1}{2}$. The product $q\|q\varphi\|$ appears to converge toward $\frac{1}{\sqrt{5}} \approx 0.44721359$. For the odd approximations (which are overestimates), the product is strictly increasing, while for the even approximations (which are underestimates), the product is strictly decreasing. They split up accordingly in table 6.

These results allow us to glimpse into Hurwitz' theorem. Dirichlet's theorem states that for an irrational number $\alpha$, there are infinitely many integers $q$ such that $q\|q\alpha\| < 1$. For the Golden ratio, there are infinitely many integers $q$ such that $q\|q\varphi\| < \frac{1}{\sqrt{5}}$, and this constant cannot be lowered. As all the partial quotients of the Golden ratio are equal to 1, we can surmise that this number is of all irrational numbers the most badly approximable, and that for an irrational number $\alpha$, there are infinitely many integers $q$ such that $q\|q\alpha\| < \frac{1}{\sqrt{5}}$, where the constant $\frac{1}{\sqrt{5}}$ cannot be improved. This is Hurwitz' theorem, which in fact can be proved without the use of continued fraction theory.

Continued fraction theory and *BRAIN* work jointly to provide a unified picture of the best rational approximation problem. Given an irrational number $\alpha$ and integer $q \neq 0$ the relation $\left| \alpha - \frac{p}{q} \right| < f(q)$ possesses a unique solution $p = [q\alpha]$ when $f(q) < \frac{1}{2q}$, and ceases to have an infinite number of solutions when $f(q) < \frac{1}{\sqrt{5}q^2}$. The *BRAIN* theorems and algorithms, together with the theorems and algorithms of continued fraction theory serve to establish a bridge between these two extremes. For instance, if $\left| \alpha - \frac{p}{q} \right| < \frac{1}{2q^2}$, then $\frac{p}{q}$ is a convergent (this explains the identity of type I and type II best approximations for the Golden ratio). In terms of algorithms, the RCF algorithm can be enriched with the NICF or Minkoswki's diagonal continued fraction algorithms; the algorithm for type I best approximations is completed by the algorithm of continued fraction theory derived by Lagrange, Stephen Smith and their followers.

## 6. Concluding remarks

There are several possible introductions to the field of Diophantine approximation: one can start with continued fractions; another approach, common is the literature on the geometry of numbers, proceeds from Dirichlet's theorem; yet another possibility is to use Farey sequences. The present work offers a complementary approach, simpler than that of the continued fraction, as attested by the rudimentary character of the algorithms described here when compared to the ingenuity of the algorithms of Euclid and of the continuous fraction apparatus. Its central results, *BRAIN II* theorem and its dedicated algorithm, are readily derived from the *TRAIN* lemma. Somewhat unexpectedly, these results are not found in the bibliography given below. This raises stimulating questions for a historian of ideas.

## Acknowledgment

This work has been carried out for the introductory course of the MD-PhD curriculum of the Université Pierre et Marie Curie. I thank Florent Soubrier, creator of the curriculum in 2008 and its manager until last May for his invitation to participate in this project, his initial trust and constant support over the last ten years.

Khinchin Continued fractions 1963; Dirichlet principle in the theory of Diophantine approximations (in Russian) 1948

Khovanskii The application of continued fractions to problems in approximation theory 1963

Knott Ron Knott's webpages on Mathematics, various articles (Introduction to Continued Fractions, Fibonacci Numbers and the Golden Section). June 2018

Knuth The Art of Computer Programming Fundamental Algorithms Vol. 1 1973; The distribution of continued fractions approximations 1984

Knuth & Patashnik Concrete Mathematics 1994

Koblitz - A course in number theory and cryptography 1994

Koksma Diophantische Approximationen 1936

Koshy Fibonacci and Lucas numbers with applications 2001

Kraaikamp & Liardet Good Approximations and Continued Fractions 1991

Krushchev The great theorem of A.A. Markoff and Jean Bernoulli sequences 2010

Kuipers & Niederreiter Uniform distribution of sequences 1974

Lagrange Oeuvres Tomes II 1868 et VII 1877

Lang Report on Diophantine approximations 1965

Lang Introduction to Diophantine Approximations 1995

LeVeque Topics In Number Theory Volume 1 1955

Lucas Théorie des nombres 1891

Manin & Panchishkin Introduction to Modern Number Theory Fundamental Problems, Ideas and Theories 2005

Markoff Sur une question de Jean Bernoulli 1882

MacTutor History of Mathematics archive website. Various articles. June 2018

Merton On the Shoulders of Giants 1965

Milnor Dynamics in One Complex Variable 1999

Minkowski Ueber Eigenschaften von ganzen Zahlen, die durch räumliche Anschauung erschlossen sind 1893; Diophantische Approximationen 1907; Geometrie der Zahlen 1910 and other texts

Mordell Diophantine equations 1969

Moore (Charles) An introduction to continued fractions 1964

Moore (Sam) On Rational Approximation 2013

Moser On commuting circle mappings and simultaneous Diophantine approximations 1990

Moshchevitin Continued fraction 1999; On some open problems in Diophantine approximation 2012

Narkiewicz Rational number theory in the 20th Century From PNT to FLT 2012

Neuenschwander Cryptology 2004

Nikishin and Sorokin - Rational Approximations and Orthogonality 1991

Niven Irrational numbers 1956; Diophantine approximations 1963

Niven Zuckerman & Montgomery An Introduction to the Theory of Numbers 1991

Ofman Une leçon de mathématiques en 399 BCE (sd)

Olds Continued fractions 1963

Olds Lax & Davidoff Geometry of Numbers 2000

**Tables**

| $q$ | $[q\pi]$ | $\|q\pi\|/q$ |
|---|---|---|
| 339 | 1065- | 2.66764E-07 |
| 678 | 2130- | 2.66764E-07 |
| 791 | 2485- | 2.66764E-07 |
| 113 | 355- | 2.66764E-07 |
| 226 | 710- | 2.66764E-07 |
| 452 | 1420- | 2.66764E-07 |
| 904 | 2840- | 2.66764E-07 |
| 565 | 1775- | 2.66764E-07 |
| 897 | 2818+ | 9.59896E-06 |
| 911 | 2862- | 9.98088E-06 |
| 784 | 2463+ | 1.10209E-05 |
| 798 | 2507- | 1.13564E-05 |
| 671 | 2108+ | 1.29218E-05 |
| 685 | 2152- | 1.31858E-05 |
| 558 | 1753+ | 1.55927E-05 |
| 572 | 1797- | 1.5738E-05 |
| 459 | 1442- | 1.95468E-05 |
| 918 | 2884- | 1.95468E-05 |
| 445 | 1398+ | 1.96199E-05 |
| 890 | 2796+ | 1.96199E-05 |

Table 1. Computation of $BRAIN\ I(\pi, k_0)$ for the first $1000$ integers (corresponding to a denominator $1 \leq q \leq 1000$) with a spreadsheet made of 3 columns of $1000$ cells. Each row contains $q$ in the first column, $[q\pi]$ in the second and $\|q\pi\|/q$ in the third. The spreadsheet has been sorted using the cells of the third column, ordered from the smallest to the largest. Only the 20 best results are displayed. The symbols plus and minus in the second row are used to indicate fractions which are underestimates $(\pi - \frac{[q\pi]}{q}) > 0$ and overestimates $(\pi - \frac{[q\pi]}{q}) < 0$ of $\pi$ respectively. The five fractions $\frac{1398}{445}, \frac{1753}{558}, \frac{2108}{671}, \frac{2463}{784}, \frac{2818}{897}$ correspond to best underestimates of the first kind of $\pi$: they minimize $\min_{1 \leq q \leq q_k} \frac{[q\pi]}{q}$ instead of $\min_{1 \leq q \leq q_k} \frac{\|q\pi\|}{q}$.

|  | $[q\pi]$ | $\|q\pi\|$ |
|---|---|---|
| 113 | 355- | 3.01444E-05 |
| 226 | 710- | 6.02887E-05 |
| 339 | 1065- | 9.04331E-05 |
| 452 | 1420- | 0.000120577 |
| 565 | 1775- | 0.000150722 |
| 678 | 2130- | 0.000180866 |
| 791 | 2485- | 0.00021101 |
| 904 | 2840- | 0.000241155 |
| 897 | 2818+ | 0.00861027 |
| 784 | 2463+ | 0.008640414 |
| 671 | 2108+ | 0.008670559 |
| 558 | 1753+ | 0.008700703 |
| 445 | 1398+ | 0.008730847 |
| 332 | 1043+ | 0.008760992 |
| 219 | 688+ | 0.008791136 |
| 106 | 333+ | 0.008821281 |
| 7 | 22- | 0.008851425 |
| 120 | 377- | 0.008881569 |
| 233 | 732- | 0.008911714 |
| 346 | 1087- | 0.008941858 |

Table 2. Computation of $BRAIN\ II(\pi, k_0)$ for the first $1000$ integers (corresponding to a denominator $1 \le q \le 1000$) with a spreadsheet made of 3 columns of $1000$ cells. Each row contains $q$ in the first column, $[q\pi]$ in the second and $\|q\pi\|$ in the third. The spreadsheet has been sorted using the cells of the third column, ordered from the smallest to the largest. Only the 20 best results are displayed. The symbols plus and minus in the second row are used to indicate fractions which are underestimates $(\pi - \frac{[q\pi]}{q}) > 0$ and overestimates $(\pi - \frac{[q\pi]}{q}) < 0$ of $\pi$ respectively.

| $q$ | $[q\pi]$ | $q\|q\pi\|$ |
|---|---|---|
| 113 | 355- | 0.003406312 |
| 226 | 710- | 0.013625248 |
| 339 | 1065- | 0.030656808 |
| 452 | 1420- | 0.054500992 |
| 7 | 22- | 0.061959974 |
| 565 | 1775- | 0.085157799 |
| 678 | 2130- | 0.122627231 |
| 1 | 3+ | 0.141592654 |
| 791 | 2485- | 0.166909287 |
| 904 | 2840- | 0.218003966 |
| 14 | 44- | 0.247839896 |
| 21 | 66- | 0.557639767 |
| 2 | 6+ | 0.566370614 |
| 6 | 19- | 0.902664471 |
| 106 | 333+ | 0.935055735 |
| 28 | 88- | 0.991359586 |

Table 3. Computation of $BRAIN\ III(\pi, k_0)$ for the first $1000$ integers (corresponding to a denominator $1 \leq q \leq 1000$) with a spreadsheet made of 3 columns of $1000$ cells. Each row contains $q$ in the first column, $[q\pi]$ in the second and $q\|q\pi\|$ in the third. The spreadsheet has been sorted using the cells of the third column, ordered from the smallest to the largest. All the results corresponding to $q\|q\pi\| < 1$ are displayed. The symbols plus and minus in the second row are used to indicate fractions which are underestimates $(\pi - \frac{[q\pi]}{q}) > 0$ and overestimates $(\pi - \frac{[q\pi]}{q}) < 0$ of $\pi$ respectively.

| $q$ | $[q\varphi]$ | $\|q\varphi\|/q$ |
|---|---|---|
| 987 | 1597- | 4.59072E-07 |
| 610 | 987+ | 1.20186E-06 |
| 843 | 1364+ | 3.14652E-06 |
| 377 | 610- | 3.14653E-06 |
| 754 | 1220- | 3.14653E-06 |
| 898 | 1453- | 6.10034E-06 |
| 233 | 377+ | 8.23768E-06 |
| 466 | 754+ | 8.23768E-06 |
| 932 | 1508+ | 8.23768E-06 |
| 699 | 1131+ | 8.23768E-06 |
| 521 | 843- | 8.23774E-06 |
| 665 | 1076- | 1.1124E-05 |
| 809 | 1309- | 1.29828E-05 |
| 788 | 1275+ | 1.36842E-05 |
| 953 | 1542- | 1.42799E-05 |
| 555 | 898+ | 1.59707E-05 |
| 877 | 1419+ | 1.80252E-05 |
| 966 | 1563+ | 2.15664E-05 |
| 322 | 521+ | 2.15664E-05 |
| 644 | 1042+ | 2.15664E-05 |

Table 4. Computation of $BRAIN\ I(\varphi, k_0)$ for the first $1000$ integers (corresponding to a denominator $1 \leq q \leq 1000$) with a spreadsheet made of 3 columns of $1000$ cells. Each row contains $q$ in the first column, $[\varphi\pi]$ in the second and $\|\varphi\pi\|/\varphi$ in the third. The spreadsheet has been sorted using the cells of the third column, ordered from the smallest to the largest. Only the 20 best results are displayed. The symbols plus and minus in the second row are used to indicate fractions which are underestimates $(\varphi - \frac{[q\varphi]}{q}) > 0$ and overestimates $(\pi - \frac{[q\varphi]}{q}) < 0$ of $\varphi$ respectively. Four best approximations can be determined in the first twenty best results: $\frac{377^+}{233}, \frac{610^-}{377}, \frac{987^+}{610}, \frac{1597^-}{987}$.

| $q$ | $[q\varphi]$ | $\|q\varphi\|$ |
|---|---|---|
| 987 | 1597- | 0.000453104 |
| 610 | 987+ | 0.000733137 |
| 377 | 610- | 0.001186241 |
| 233 | 377+ | 0.001919379 |
| 754 | 1220- | 0.002372483 |
| 843 | 1364+ | 0.002652516 |
| 144 | 233- | 0.00310562 |
| 466 | 754+ | 0.003838757 |
| 521 | 843- | 0.004291861 |
| 89 | 144+ | 0.005024999 |
| 898 | 1453- | 0.005478103 |
| 699 | 1131+ | 0.005758136 |
| 288 | 466- | 0.00621124 |
| 322 | 521+ | 0.006944377 |
| 665 | 1076- | 0.007397481 |
| 932 | 1508+ | 0.007677515 |
| 55 | 89- | 0.008130619 |
| 555 | 898+ | 0.008863756 |
| 432 | 699- | 0.00931686 |
| 178 | 288+ | 0.010049997 |

Table 5. Computation of $BRAIN\ II(\varphi, k_0)$ for the first $1000$ integers (corresponding to a denominator $1 \leq q \leq 1000$) with a spreadsheet made of 3 columns of $1000$ cells. Each row contains $q$ in the first column, $[q\varphi]$ in the second and $\|q\varphi\|$ in the third. The spreadsheet has been sorted using the cells of the third column, ordered from the smallest to the largest. Only the 20 best results are displayed. The symbols plus and minus in the second row are used to indicate fractions which are underestimates $(\varphi - \frac{[q\varphi]}{q}) > 0$ and overestimates $(\pi - \frac{[q\varphi]}{q}) < 0$ of $\pi$ respectively. Seven best approximations can be determined in the first twenty best results: $\frac{89^-}{55}, \frac{144^+}{89}, \frac{233^-}{144}, \frac{377^+}{233}, \frac{610^-}{377}, \frac{987^+}{610}, \frac{1597^-}{987}$.

| $q$ | $[q\varphi]$ | $q\|q\varphi\|$ |
|---|---|---|
| 1 | 2- | 0.381966011 |
| 3 | 5- | 0.437694101 |
| 8 | 13- | 0.44582472 |
| 21 | 34- | 0.447010961 |
| 55 | 89- | 0.447184032 |
| 144 | 233- | 0.447209282 |
| 377 | 610- | 0.447212966 |
| 987 | 1597- | 0.447213504 |
| 610 | 987+ | 0.447213836 |
| 233 | 377+ | 0.447215243 |
| 89 | 144+ | 0.447224888 |
| 34 | 55+ | 0.447290995 |
| 13 | 21+ | 0.447744099 |
| 5 | 8+ | 0.450849719 |
| 2 | 3+ | 0.472135955 |

Table 6. Computation of $BRAIN\ III(\varphi, k_0)$ for the first $1000$ integers (corresponding to a denominator $1 \le q \le 1000$) with a spreadsheet made of 3 columns of $1000$ cells. Each row contains $q$ in the first column, $[\varphi\pi]$ in the second and $\varphi\|q\varphi\|$ in the third. The spreadsheet has been sorted using the cells of the third column, ordered from the smallest to the largest. All the results corresponding to $\varphi\|q\varphi\| < 1$ are displayed. The symbols plus and minus in the second row are used to indicate fractions which are underestimates $(\varphi - \frac{[q\varphi]}{q}) > 0$ and overestimates $(\varphi - \frac{[q\varphi]}{q}) < 0$ of $\pi$ respectively.